\documentclass[12pt]{article}

\usepackage{amsmath}
\usepackage{amssymb}
\usepackage[mathscr]{eucal}
\usepackage{graphicx}
\usepackage{color}

\usepackage{caption}
\usepackage{subcaption}

\newtheorem{Theorem}{\sc Theorem}

\newtheorem{Lemma}[Theorem]{\sc Lemma}

\newtheorem{Remark}[Theorem]{\sc Remark}
\newtheorem{Example}[Theorem]{\sc Example}

\newcommand{\cQ}{\mbox{{${\cal Q}$}}}

\newcommand{\cP}{\mbox{{${\cal P}$}}}

\newcommand{\ve}{\mbox{{$\varepsilon$}}}
\newcommand{\R}{\mbox{{${\mathbb{R}}$}}}

\newcommand{\cF}{\mbox{{${\cal F}$}}}

\newcommand{\bu}{\mbox{\boldmath{$u$}}}
\newcommand{\bv}{\mbox{\boldmath{$v$}}}

\newcommand{\bx}{\mbox{\boldmath{$x$}}}

\newcommand{\fb}{\mbox{\boldmath{$f$}}}

\newcommand{\bxi}{\mbox{\boldmath{$\xi$}}}

\newcommand{\bsigma}{\mbox{\boldmath{$\sigma$}}}

\newcommand{\btau}{\mbox{\boldmath{$\tau$}}}
\newcommand{\bvarepsilon}{\mbox{\boldmath{$\varepsilon$}}}

\newcommand{\bnu}{\mbox{\boldmath{$\nu$}}}

\newcommand{\bzero}{\mbox{\boldmath{$0$}}}
\newcommand{\la}{\langle}
\newcommand{\ra}{\rangle}

\def\sqr#1#2{{
    \vcenter{
         \vbox{\hrule height.#2pt
               \hbox{\vrule width.#2pt height#1pt \kern#1pt
                     \vrule width.#2pt
               }
               \hrule height.#2pt
         }
    }
}}

\def\real{\mathbb{R}}

\def\lista#1
{{ \itemindent 0.0cm \labelsep .2cm \leftmargin 0.8cm \rightmargin
0.0cm \labelwidth 0.6cm \topsep 0.0mm
\parsep 0.0mm
\itemsep 0.0mm
\begin{list}{}
{ \setlength{\leftmargin}{.8cm} \setlength{\rightmargin}{0.0cm}
\setlength{\parsep}{0.0mm} \setlength{\topsep}{.0mm}
\setlength{\parskip}{.0cm} \setlength{\itemsep}{.0cm} }
{#1}\end{list}} }

\newcounter{theorem}

\begin{document}

\title{\bf  A new penalty method for elliptic quasivariational inequalities}

\vspace{22mm}
{\author{Piotr Bartman-Szwarc$^{1, 2}$,\ Anna Ochal$^{2}$\footnote{Corresponding author, E-mail : anna.ochal@uj.edu.pl },\ Mircea Sofonea$^{3}$\\ and\\ Domingo A. Tarzia$^{4,5}$\\[6mm]	
{\it \small $1$ Doctoral School of Exact and Natural Sciences, Jagiellonian University in Krakow}
                \\{\it\small Lojasiewicza 11, 30348 Krakow, Poland}
                \\ [4mm]
                {\it \small $2$ Chair of Optimization and Control, Jagiellonian University in Krakow}\\
{\it\small Lojasiewicza 6, 30348 Krakow, Poland}		\\[6mm]	
{\it \small $^3$ Laboratoire de Math\'ematiques et Physique, University of Perpignan Via Domitia}
\\{\it\small 52 Avenue Paul Alduy, 66860 Perpignan, France}		\\[6mm]		
{\it\small $^{4}$  Departamento de Matem\'atica, FCE, Universidad Austral}\\
		{\it \small Paraguay 1950, S2000FZF Rosario, Argentina}\\
{\it\small $^{5}$ CONICET, Argentina}}}

\date{}
\maketitle
\thispagestyle{empty}

\vskip 5mm

\noindent {\small{\bf Abstract.}
 We consider a class of elliptic quasivariational  inequalities in a reflexive Banach space~$X$ for which we recall a convergence criterion obtained in \cite{GOST}. Each inequality $\cP$ in the class is governed by a set of constraints $K$ and has a unique solution $u\in K$. The criterion  provides necessary and sufficient conditions which guarantee that an arbitrary sequence $\{u_n\}\subset X$ converges to the solution $u$.
Then, we consider a sequence $\{\cP_n\}$ of unconstrained variational-hemivariational inequalities  governed by a sequence of parameters~$\{\lambda_n\}\subset\real_+$.
We use our criterion to deduce that, if for each $n\in\mathbb{N}$  the term $u_n$ represents a solution of Problem $\cP_n$, then  the sequence
$\{u_n\}$ converges to $u$ as $\lambda_n\to 0$.
We apply our abstract results in the study of an elastic frictional contact problem with unilateral constraints and provide the corresponding mechanical interpretations. We also present numerical simulation in the study of a two-dimensional example which represents an evidence of our convergence results.}


\vskip 4mm


\section{Introduction}\label{s1}
\setcounter{equation}0


\medskip\noindent
A large number of mathematical models in Physics, Mechanics and Engineering Science lead to nonlinear problems  governed by unilateral constraints.  The famous Signorini contact problem which describes the equilibrium of an elastic body in contact with a rigid foundation and the heat transfer problem across a semipermeable membrane are two examples among others.
The penalty method represents one of the most popular methods in the study of constrained problems.
Its main ingredients are the following:
one  replaces the statement of the constrained problem with a sequence of unconstrained problems, the so-called penalty problems, which are governed by a sequence of  penalty parameters, say $\{\lambda_n\}$;  one proves that each unconstrained problem has at least one  solution $u_n$; finally,  one proves that the sequence $\{u_n\}$ converges to the unique solution $u$ of the constrained problem as $\{\lambda_n\}$ converges.

The interest in penalty methods is two folds. First, solving an unconstrained problem is more convenient from numerical point of view. Second, from theoretical point of view the penalty methods establish the link between problems which could have a different feature. For instance, the penalty method allows us to approach the solution of a  contact problem with a rigid foundation by the solution of a contact problem with a deformable foundation, for a small deformability coefficient.

\medskip
In this framework very general, the following two questions arise.

\medskip
a) Given Problem $\cP$ with unilateral constraints which has a unique solution $u$, how to choose the sequence of the penalty Problems $\{\cP_n\}$  governed by a sequence of penalty parameters $\{\lambda_n\}$?

\medskip
b) How to guarantee that a solution $u_n$ of Problem $\cP_n$ converges to the solution $u$ of Problem $\cP$?

\medskip
The answer to these questions varies from problem to problem, depends on the assumption made on the data and requires on good knowledge of the convergence results to the solution of the original problem.
Nevertheless, some elementary examples show that concerning question a) the choice is not unique, i.e., it is possible to construct several penalty problems for the same constrained problem. Concerning question b) we stress that it would be interesting to describe a general framework in which several convergence results for penalty method can be cast.

\medskip
In this paper we provide an answer to the questions above, in the case when Problem~$\cP$ is an elliptic  quasivariational inequality. The functional framework is the following:   $X$~is a reflexive Banach space endowed with the norm  $\|\cdot\|_X$, $X^*$ represents its dual space and  $\la\cdot,\cdot\ra$ denotes the duality pairing mapping; also,
$K\subset X$, $A\colon X\to X^*$, $\varphi\colon X\to\R$ and $f\in X^*$. Then, the inequality problem we consider in this paper is as follows.

\medskip\medskip\noindent{\bf Problem}  ${\cal P}$. {\it Find $u$ such that}
\begin{equation}\label{1}u\in K,\qquad\la Au,v-u\ra+\varphi(v)-\varphi(u) \ge\la f,v-u\ra \qquad\forall\,v\in K.
\end{equation}

Elliptic quasivariational inequalities of the above form have been studied by many authors, under different assumptions on the data.
Existence and uniqueness results can be found in \cite{BC, B, G, Kind-St}, for instance.
The proofs provided there are based on
arguments of monotonicity and convexity, including properties of the subdifferential
of a convex function.   References on the approximation of Problem $\cP$, including penalty methods and error estimates for discrete finite element schemes can be found in \cite{HS}.
Applications of variational inequalities of the form (\ref{1}) in Mechanics can be found in
the books \cite{C,DL,EJK,HS,HHNL,KO, Pa} and, more recently, \cite{S}.

In this paper we describe the convergence of a sequence $\{u_n\}\subset X$ to the solution $u$ of the variational inequality (\ref{1}). More precisely, we recall  a {\it convergence criterion} provided in \cite{GOST}, which states necessary and sufficient conditions for the convergence $u_n\to u$ in $X$. Our aim in this paper is two folds.  The first one is to provide a penalty method in the study of \eqref{1} based on the study of a class of unconstrained variational-hemivariational inequalities.   Variational-hemivariational inequalities  represent a special class of inequality problems  which have both a convex and nonconvex structure. They have been intensively studied in the last decades as shown in \cite{HS-Acta, MOS30, NP, Pana1, P, PK, SMBOOK} and the references therein. The penalty method we introduce in the current paper is new since, at the best of our knowledge, the penalty problems of \eqref{1} studied in the literature are in a form of quasivariational inequality governed by a penalty operator and the main ingredients used there are based on pseudomonotonicity arguments. In contrast, in our approach we use as a penalty term the directional derivative of a nonsmooth function. Even if the penalty method we introduce here could be unconvenient from numerical point of view, it has the merit to represent an alternative which, through the example of   inequality \eqref{1}, illustrates that penalty methods are not unique.
Finally, our second aim is to show how our theoretical results can be applied in the study of mathematical models of contact which, in a~variational formulation, leads to a quasivariational inequality of the form \eqref{1}.  The mechanical interpretations related to these applications, supported by numerical simulations, fully justify the interest in the  penalty method we introduce in this work.

The rest of the manuscript is structured as follows. In Section \ref{s2} we present the preliminary material we need in the rest of the paper. In particular, we prove the existence result for variational-hemivariational inequality with a parameter $\lambda>0$ and recall a criterion of convergence to the solution $u$ of inequality (\ref{1}), Theorem \ref{tm}.
In Section \ref{s4} we introduce a~sequence of unconstrained variational-hemivariational inequalities, governed by a penalty parameter. We prove the solvability of each inequality and the convergence of any sequence of solutions to the solution of \eqref{1}, as the penalty parameter converges to zero. An application in Contact Mechanics is provided in Section  \ref{s5}, together with the corresponding mechanical interpretations.
Finally, in Section~\ref{s6} we present numerical simulations which validate our theoretical convergence results.

\section{Preliminaries}\label{s2}
\setcounter{equation}0

Everywhere in this paper, unless it is specified otherwise, we use the functional framework described in Introduction. Notation $0_X$ represents the zero element of $X$. All the limits, upper and lower limits below are considered as $n\to\infty$, even if we do not mention it explicitly. The symbols ``$\rightharpoonup$"  and ``$\to$"
denote the weak and the strong convergence in various spaces which will be specified, except in the case when these convergences take place in $\real$.  For a sequence $\{\ve_n\}\subset\R_+$ which  converges to zero we use the short hand notation $0\le\ve_n\to 0$. Finally,
we denote by $d(u,K)$ the distance between an element $u\in X$ to the set $K$, that is
\begin{equation}\label{d}
	d(u,K)=\inf_{v\in K} \|u-v\|_X.
\end{equation}

\medskip
In the study of (\ref{1}) we consider the following assumptions.
\begin{eqnarray}
	&&\hspace{0mm}\label{K}
	K \ \mbox{is a nonempty closed convex subset of} \ X.
	\\[3mm]
	&&\hspace{-6mm}\label{A}
	\left\{\begin{array}{l}A \colon X \to X^* \ \mbox{is  pseudomonotone and strongly monotone, i.e.:}\\ [1mm]
		{\rm (a)} \  A\  \mbox{is bounded and}\ u_n\rightharpoonup u \ \mbox{in}\ X\\
		\qquad\mbox{with}\displaystyle \limsup\,\langle A u_n, u_n -u \rangle \le 0\ \\
		\qquad\mbox{imply}\ \ \displaystyle \liminf\, \langle A u_n, u_n - v \rangle\ge \langle A u, u - v \rangle\ ~\forall\, v \in X.\\ [1mm]
		{\rm (b)} \ \mbox{there exists}\  m_A > 0\ \mbox{such that}\\
		\qquad\langle Au - Av, u-v \rangle
		\ge
		m_A \| u - v \|_X^{2}\ \ \ \forall\,
		u,\, v \in X.
	\end{array}\right.
\\[3mm]
	&&\hspace{0mm}\label{fi}
		\varphi \colon X \to \real \ \mbox{is convex and lower semicontinuous.}
	\\[3mm]
	&&\hspace{0mm}\label{f}
	f \in X^*.
\end{eqnarray}

\medskip

\medskip

The unique solvability of inequality (\ref{1}) is given by a general result, proved in \cite{SMBOOK}, for instance.

\begin{Theorem}\label{t1}
	Assume  $(\ref{K})$--$(\ref{f})$. Then, inequality~$(\ref{1})$ has a unique solution $u \in K$.
\end{Theorem}

We now move to a different type of inequalities, the so-called variational-hemivariational inequalities.  To this end, besides the data already used in the statement of Problem $\cP$, we consider a reflexive space $Y$, a locally Lipschitz continuous function $j\colon Y\to \real$ and an operator $\gamma\colon X\to Y$. Then, denoting by $j^0(u;v)$ the Clarke directional derivative of $j$ at the point $u$ in the direction $v$, the inequality problem we consider with a fixed $\lambda >0$, is the following.

\medskip\medskip\noindent{\bf Problem}  ${\cal Q}$. {\it Find $u$ such that}
\begin{equation}\label{2}u\in X,\quad\la Au,v-u\ra+\varphi(v)-\varphi(u)+\frac{1}{\lambda}\,j^0(\gamma u; \gamma v-\gamma u) \ge\la f,v-u\ra \quad\forall\,v\in X.
\end{equation}

We recall shortly that for a locally Lipschitz function $j\colon Y\to \real$ we define (see \cite{Cl}) the generalized (Clarke) directional derivative of $j$ at the point $u \in Y$ in the direction $v \in Y$ by
\begin{align*}
  &j^0(u;v) := \limsup_{w \to u, \lambda \searrow 0} \frac{j(w+\lambda v)  - j(w)}{\lambda}.
\end{align*}
The generalized gradient of $j$ at $u$ is a subset of the dual space $Y^*$ given by
\begin{align*}
  &\partial j(u) := \{\xi \in Y^*\, | \, \langle \xi, v\rangle_{Y^*\times Y}\leq j^0(u;v) \  \mbox{ for all } v \in Y \}
\end{align*}
and $\partial j:Y\to 2^{Y^*}$ represents the Clarke subdifferential of the function $j$.
The function $j$ is said to be regular (in the sense
of Clarke)  if for all $u,\, v \in Y$ the one-sided directional derivative
\begin{equation*}
	j'(u; v) = \lim_{\lambda \downarrow 0}
	\frac{j(u + \lambda v) - j(u)}{\lambda}
\end{equation*}
exists and $j^0(u; v) = j'(u; v)$.

\medskip

Variational-hemivariational inequalities of the form \eqref{2} have been intensively studied in the literature. Reference in the field include \cite{Han, LM, MNZ, MYZ, MOS30, SMBOOK}.  There, various existence, uniqueness and convergence results  have been obtained, under various assumptions on the data.
Here, in the study of  (\ref{2}) we assume the following.

\begin{eqnarray}
	&&\hspace{-6mm}\label{j} \left\{
	\begin{array}{ll}
		j \colon Y \to \real \ \mbox{is such that:}\\ [1mm]
		{\rm (a)} \
		j \ \mbox{is locally Lipschitz.} \\ [1mm]
		{\rm (b)} \
		\| \xi \|_{Y^*} \le c_0 + c_1 \, \| v \|_Y \ \ \mbox{for all}
		\ v \in Y,\ \xi\in \partial j(v)~\\
		\qquad \mbox{with} \ c_0, c_1 \ge 0. \\ [1mm]
		{\rm (c)} \
		\mbox{there exists} \ d_1 \ge 0 \ \mbox{such that} \\ [1mm]
		\qquad
        j^0(v;  - v)  \le d_1\, (1+\|v\|_Y) \  \ \ \forall\, v \in Y.
	\end{array}
	\right.
\\[3mm]
	&&\hspace{-6mm}\label{ga}
	\mbox{\quad\ $\gamma\colon X\to Y$ is a linear compact operator.}
	\end{eqnarray}

\medskip
Under these assumptions we have the following existence result.

\begin{Theorem}\label{t2}
	Assume  $(\ref{A})$--$(\ref{f})$, $(\ref{j})$ and $(\ref{ga})$. Then, inequality~$(\ref{2})$ has at least one solution $u \in X$.
\end{Theorem}

\noindent {\it Proof.}
We consider the following inclusion: find $u\in X$ such that
\begin{align}\label{q1}
Au + \frac{1}{\lambda}\, \gamma^* \partial j(\gamma u) + \partial \varphi (u) \ni f,
\end{align}
where $\gamma^*\colon Y^*\to X^*$ is the adjoint operator to $\gamma$ and $\partial \varphi$ represents the convex subdifferential of $\varphi$.
We introduce two multivalued operators
\begin{align*}
  T_1\colon X\to 2^{X^*}, &  \qquad T_1(v)=Av+\frac{1}{\lambda}\,\gamma^*\partial j(\gamma v) \quad \mbox{for\ } v\in X, \\
  T_2\colon X\to 2^{X^*}, & \qquad T_2(v)=\partial \varphi (v) \quad \mbox{for\ } v\in X,
\end{align*}
respectively. It is well known that the operator $T_2$ is maximal monotone on $X$.

We prove that the operator $T_1$ is bounded, coercive and pseudomonotone (\cite[Definition~3.57]{MOSBOOK}). To this end, we consider the operator
$B \colon  X \to 2^{X^*}$
given by
\begin{align*}
Bv = \frac{1}{\lambda}\, \gamma^* \partial j (\gamma v) \ \ \mbox{for all} \ \ v \in X.
\end{align*}
We show that $B$ is pseudomonotone and bounded.
The boundedness of $B$ follows easily from \eqref{j}(b)
\begin{align*}
\|Bv\|_{X^*}\leq \frac{1}{\lambda}\, \|\gamma^*\|\,\|\partial j(\gamma v)\|_{Y^*}\leq \frac{1}{\lambda}\, \|\gamma^*\|\,(c_0+c_1 \|\gamma\|\,\| v\|_{X})
\ \ \mbox{for \ all \ } v\in X,
\end{align*}
where $\|\gamma \|$ denotes the norm of linear operator $\gamma$. Moreover, we observe that
the values of $\partial j$ are nonempty, convex and weakly compact subsets of $Y^*$, \cite[Proposition~3.23(iv)]{MOSBOOK}.
Hence, for all $v\in X$ the set $Bv$ is nonempty, closed and convex in $X^*$. Next, using \cite[Proposition~3.58(ii)]{MOSBOOK}, it is enough to prove that $B$ is generalized pseudomonotone in order to show that $B$ is pseudomonotone.

Let $v_n$, $v \in X$, $v_n \rightharpoonup v$ in $X$, $v^*_n$, $v^* \in X^*$, $v^*_n \rightharpoonup v^*$ in $X^*$,
$v_n^* \in Bv_n$ and\
$\limsup\, \langle v_n^*, v_n - v \rangle \le 0$.
We show that
$$
v^* \in Bv \ \ \ \mbox{and} \ \ \
\langle v_n^*, v_n \rangle \to \langle v^*, v \rangle.
$$
Since $v_n^* \in Bv_n$, we have
\begin{align}\label{q2}
  v_n^* = \frac{1}{\lambda}\, \gamma^* \zeta_n & \ \ \mbox{with \ } \zeta_n \in \partial j(\gamma v_n).
\end{align}
And, from \eqref{j}(b), it follows that $\{\zeta_n\}$ is bounded in $Y^*$. Hence,
we may assume that, at least for a subsequence, denoted in the same way,
we have
\begin{align}\label{q3}
  \zeta_n \rightharpoonup \zeta & \ \ \mbox{in \ }Y^*.
\end{align}
Moreover, by the compactness of the operator~$\gamma$, we obtain
\begin{align}\label{q4}
  \gamma v_n \to \gamma v & \ \ \mbox{in \ } Y.
\end{align}
Hence, from \eqref{q2}-\eqref{q4}, it is obvious that
\begin{align*}
	\langle v_n^*, v_n \rangle =
	\frac{1}{\lambda}
	\langle \zeta_n, \gamma v_n \rangle_{Y^*\times Y}
	\longrightarrow \frac{1}{\lambda}
	\langle \zeta, \gamma v \rangle_{Y^*\times Y} =
	\langle \frac{1}{\lambda} \gamma^* \zeta, v \rangle =
	\langle v^*, v \rangle .
\end{align*}

Next, using the strong-weak closedness of the graph of $\partial j$,
(\cite[Proposition~3.23]{MOSBOOK}), we get $\zeta \in \partial j(\gamma v)$.
Finally, by $v_n^* = \frac{1}{\lambda}\, \gamma^* \zeta_n$, we have $v^* = \frac{1}{\lambda}\, \gamma^* \zeta$,
and consequently
\begin{align*}
  v^* \in \frac{1}{\lambda}\, \gamma^* \partial j(\gamma v) = Bv.
\end{align*}
This shows that $B$ is generalized pseudomonotone and also $B$ is pseudomonotone. Moreover, since
$A \colon X \to X^*$ is pseudomonotone (and therefore, it is pseudomonotone and bounded as a multivalued operator), using \cite[Proposition~3.59(ii)]{MOSBOOK}, we obtain that $T_1=A+B$ is bounded and pseudomonotone.

Next, we claim that the operator $T_1$ is coercive. To this end, we observe that since $A$ is strongly monotone, we have for all $v\in X$
\begin{align}\label{q5}
 \la Av, v\ra = \la Av- A0_X,v\ra+\la A0_X, v\ra \geq m_A \|v\|_X^2-\|A0_X\|_{X^*}\|v\|_X.
\end{align}
Now, let $v \in X$, $v^* \in Bv$. Thus,
$v^* = \frac{1}{\lambda}\, \gamma^* \zeta$ with $\zeta \in \partial j(\gamma v)$.
By the definition of the generalized gradient and \eqref{j}(c), we get
\begin{align}\label{q6}
\frac{1}{\lambda}\, \langle \zeta, -\gamma v \rangle_{Y^*\times Y}
\le
\frac{1}{\lambda}\,  j^0(\gamma v; - \gamma v) \le \frac{1}{\lambda}\, d_1 (1+\|\gamma\|\,\| v \|_{X}).
\end{align}
Hence, from \eqref{q5} and \eqref{q6}, we obtain
\begin{align*}
\la T_1v,v\ra=\la Av,v\ra+\la Bv,v\ra \geq m_A \|v\|_X^2-\|A0_X\|_{X^*}\|v\|_X - \frac{d_1}{\lambda}\, \|\gamma\|\,\| v \|_{X}- \frac{d_1}{\lambda},
\end{align*}
which proves that $T_1$ is coercive.

To conclude, since $T_2$ is maximal monotone on $X$ and $T_1$ is bounded, coercive and pseudomonotone, applying  a surjectivity result (\cite[Theorem~2.12]{NP}), we deduce that there exists $u \in X$ a solution to the inclusion \eqref{q1}. Finally, we observe that any solution $u\in X$ to \eqref{q1} is a solution to the inequality \eqref{2}. Indeed, if $u\in X$ solves \eqref{q1}, it means that
\begin{align*}
A u + \frac{1}{\lambda}\, \gamma^* \zeta +\xi= f \ \ \mbox{with} \ \ \zeta \in \partial j(\gamma u)\ \ \mbox{and} \ \ \xi\in \partial \varphi(u),
\end{align*}
and hence
\begin{align*}
\langle Au, v -u\rangle +
\frac{1}{\lambda} \langle \zeta, \gamma v -\gamma u\rangle_{Y^*\times Y}+ \langle \xi, v-u\rangle
= \langle f, v - u \rangle \ \ \mbox{for\ all\ } v\in X.
\end{align*}
Using the definition of subdifferentials, we have
\begin{align*}
 &\langle \zeta, \gamma v -\gamma u\rangle_{Y^*\times Y} \leq j^0(\gamma u;\gamma v-\gamma u), \\[2mm]
 &\langle \xi, v-u\rangle\leq \varphi (v)-\varphi(u),
\end{align*}
and we get
\begin{align*}
\langle Au, v -u\rangle +
\frac{1}{\lambda} j^0(\gamma u;\gamma v-\gamma u)+ \varphi (v)-\varphi(u)\leq \langle f, v - u \rangle \ \ \mbox{for\ all\ } v\in X,
\end{align*}
which means that $u\in X$ is a solution to~\eqref{2}.
This completes the proof of the theorem. \hfill$\Box$

\

We now  turn to  a convergence criterion for variational inequality $(\ref{1})$.
To this end, we need the following additional condition on the operator $A$ and function $\varphi$.
\begin{eqnarray}\label{A1}
&&\left\{\begin{array}{ll}\mbox{$A$ is a Lipschitz continuous operator, i.e., there exists $M_A>0$ such that}\\ [2mm]
	\|Au-Av\|_{X^*}\le M_A\|u-v\|_{X} \quad\forall\, u,\, v \in X.	
	\end{array}\right.	\\ [3mm]
&&\label{fi1}
	\left\{\begin{array}{ll}\mbox{For each $D>0$ there exists $L_D>0$ such that}\\ [2mm]
		|\varphi(u)-\varphi(v)|\le L_D \|u-v\|_X\quad\forall\, u,\ v,\in X\ {\rm with}\ \|u\|_X\le D, \ \|v\|_X\le D.
	\end{array}\right.	
\end{eqnarray}

\medskip
Let us notice that if \eqref{A}(b) and \eqref{A1} hold, then \eqref{A}(a) holds too.

\medskip
The following theorem introduces the convergence criterion for  $(\ref{1})$.

\begin{Theorem}\label{tm}
	Assume $(\ref{K})$--$(\ref{f})$, $(\ref{A1})$ and $(\ref{fi1})$, denote by $u$ the solution of the variational inequality $(\ref{1})$ provided by Theorem $\ref{t1}$ and let
	$\{u_n\}\subset X$. Then the following statements are equivalent:
	\begin{eqnarray}
	&&\label{c1}
	\qquad	u_n\to u\qquad{\rm in}\ X.
	\\ [5mm]
	&&\label{c2}
		\left\{\begin{array}{ll} \ {\rm (a)}\ d(u_n,K)\to 0\ ;\\[4mm]
			\ {\rm (b)}\ 	\mbox{\rm there exists $0\le\ve_n\to 0$ such that} \\[3mm]	
			\quad\quad\la Au_n,v-u_n\ra+\varphi(v)-\varphi(u_n)+\ve_n(1+\|v-u_n\|_X)\\ [1mm]
			\qquad\qquad\qquad\quad \ge\la f,v-u_n\ra  \quad\forall\,v\in K,\ n\in\mathbb{N}.
		\end{array}\right.	
	\end{eqnarray}	
\end{Theorem}

\medskip
The proof of Theorem \ref{tm} can be found in \cite[Theorem 3.1]{GOST} and is based on the following result, \cite[Lemma 3.2]{GOST}.

\begin{Lemma}\label{lm}
	Assume $(\ref{K})$--$(\ref{f})$. Then  any sequence $\{u_n\}\subset X$  which satisfies condition $(\ref{c2})\,({\rm b)}$ is bounded.
\end{Lemma}

\medskip
We remark that Theorem \ref{tm} shows that, under assumptions
$(\ref{K})$--$(\ref{f})$, $(\ref{A1})$ and
$(\ref{fi1})$, conditions $(\ref{c2})({\rm a})$ and  $(\ref{c2})({\rm b})$ represent necessary and sufficient con\-di\-tions for the convergence $(\ref{c1})$. Some elementary examples can be constructed to see that, in general, we cannot skip one of these conditions.

We end this section with the remark that Theorem \ref{tm} was obtained under the additional assumptions (\ref{A1}), (\ref{fi1})  which are not necessary in the statement of Theorem \ref{t1}. Removing or relaxing this condition is an interesting problem which clearly deserves to be investigated into future.

\section{A penalty method}\label{s4}
\setcounter{equation}0

In this section we show how Theorems \ref{t2} and \ref{tm} can be used in the study of a penalty method for inequality (\ref{1}). To this end, everywhere below  we assume that \eqref{K}--\eqref{f}, \eqref{j} and \eqref{ga} hold, even if we do not mention it explicitly.
In addition, we consider a sequence $\{\lambda_n\}\subset \real_+$ of penalty parameters.  Assuming that
$\lambda_n>0$, for each $n\in\mathbb{N}$ we consider the following unconstrained problem.

\medskip\medskip\noindent{\bf Problem}  ${\cal Q}_n$. {\it Find $u_n\in X$ such that}
\begin{equation}\label{2n}\la Au_n,v-u_n\ra+\varphi(v)-\varphi(u_n)+\frac{1}{\lambda_n}j^0(\gamma u_n; \gamma v-\gamma u_n) \ge\la f,v-u_n\ra \quad\forall\,v\in X.
\end{equation}

\noindent
Note that, under the previous assumptions, Theorem \ref{t2} guarantees that Problem~$\cQ_n$ has at least one solution. Moreover, recall that variational-hemivariational inequalities of the form
(\ref{2n}) have been considered in \cite{GMOT}, in the
study of a heat transfer problem.

\medskip
Consider now  the following additional assumptions.
\begin{eqnarray}
	&&\label{j1}j^0(\gamma u;\gamma v-\gamma u)\le 0\qquad\forall\, u\in X,\ v\in K.\\[2mm]
	&&\label{j2}u\in X,\quad j^0(\gamma u;\gamma v-\gamma u)\ge 0\quad\forall\, v\in K\quad\Longrightarrow\ u\in K.\\ [2mm]
	&&\label{lam}\lambda_n\to 0.
\end{eqnarray}

Our main result in this section is the following.

\begin{Theorem}\label{tp}
	Assume	 $(\ref{K})$--$(\ref{f})$,
	$(\ref{j})$, $(\ref{ga})$,  $(\ref{A1})$,  $(\ref{fi1})$ and $(\ref{j1})$--$(\ref{lam})$, denote by $u$ the solution of Problem $\cP$ and
	let $\{u_n\}$ be a sequence of elements in $X$ such that, for each $n\in\mathbb{N}$, $u_n$ is a solution of Problem $\cQ_n$. Then, $u_n\to u$ in $X$.
\end{Theorem}

\medskip\noindent{\it Proof.} The proof is structured in several steps, as follows.

\medskip
\noindent {\it Step {i)}  We prove that the sequence $\{u_n\}$ satisfies condition $(\ref{c2}){\rm (b)}$.}
Let $n\in\mathbb{N}$ and $v\in K$. Then, using assumption (\ref{j1}) and
inequality (\ref{2n}) we see that
\begin{equation}\label{m4}
	\la Au_n,v-u_n\ra+\varphi(v)-\varphi(u_n)\ge\la f,v-u_n\ra
\end{equation}
which shows that ${\rm (\ref{c2})(b)}$ holds with $\ve_n=0$.

\medskip
\noindent {\it Step {ii)}  We prove that any weakly convergent subsequence of the sequence $\{u_n\}$  satisfies condition ${\rm (\ref{c2})(a)}$.} Indeed,
consider a weakly convergent subsequence of the sequence $\{u_n\}$, again denoted by $\{u_n\}$.
Then, there exists an element $\widetilde{u}\in X$ such that
\begin{equation}\label{z48}
	u_n\rightharpoonup \widetilde{u}\quad{\rm in}\quad X.
\end{equation}	
We shall prove that $\widetilde{u}\in K$ and $u_n\to \widetilde{u}$ in $X$. To this end, we fix $n\in\mathbb{N}$ and $v\in X$. We use (\ref{2n}) to write	
\begin{eqnarray}
	&&\label{z50}-\,\frac{1}{\lambda_n} j^0(\gamma{u}_n;\gamma v-\gamma {u}_n)\leq \la  A{u}_n,v-{u}_n\ra+\varphi(v)-\varphi(u_n)+\la f,{u}_n-v\ra \nonumber
 \\ [2mm]
	&&\ \ \leq \|Au_n\|_{X^*}\|v-u_n\|_X+\varphi(v)-\varphi(u_n)+
	\|f\|_{X^*}\|v-u_n\|_X.
\end{eqnarray}
On the other hand, Step i) and Lemma \ref{lm} imply that the sequence $\{u_n\}$ is bounded.
Therefore, using assumption (\ref{A1}) we deduce that there exists $H>0$  such that
\begin{equation}\label{z51}
\|u_n\|_{X}\le H,\quad \|Au_n\|_{X^*}\le H\qquad\forall\, n\in\mathbb{N}.
\end{equation}
It follows now from assumption (\ref{fi1}) that
\begin{equation}\label{z52}
\varphi(v)-\varphi(u_n)\le L_H\|v-u_n\|_{X}.
\end{equation}
We now combine inequalities \eqref{z50}--\eqref{z52}
to see that
\[
-\,\frac{1}{\lambda_n} j^0(\gamma{u}_n;\gamma v- \gamma{u}_n)\le(H+L_H+\|f\|_{X^*})\|v-u_n\|_X\]
and, since, $\|v-u_n\|_X\le \|v\|_X+H$, we obtain that
\[
-j^0(\gamma{u}_n;\gamma v-\gamma{u}_n)\le\lambda_n(H+L_H+\|f\|_{X^*})(\|v\|_X+H).\]
Passing to the upper limit in the above inequality
and using assumption $\lambda_n\to 0$ we have
\begin{equation*}
\limsup\,[-j^0(\gamma{u}_n;\gamma v-\gamma{u}_n)]\leq 0
\end{equation*}
or, equivalently,
\begin{equation}\label{z53}
	0\le\liminf\,j^0(\gamma{u}_n;\gamma v-\gamma{u}_n).
\end{equation}

On the other hand, convergence \eqref{z48}, assumption \eqref{ga} and the upper semicontinuity
of the Clarke directional derivative yield
\begin{equation}\label{z54}
	\limsup\,j^0(\gamma{u}_n;\gamma v-\gamma{u}_n)\leq j^0(\gamma\widetilde{u};\gamma v-\gamma\widetilde{u}).
\end{equation}
We now combine \eqref{z53} and \eqref{z54} to see that
\begin{equation*}\label{z54n}
	0\leq j^0(\gamma\widetilde{u};\gamma v-\gamma\widetilde{u}).
\end{equation*}
Recall that this inequality holds for each $v\in K$. Therefore, using assumption \eqref{j2} we conclude that
\begin{equation}\label{z55}
\widetilde{u}\in K.
\end{equation}

Let $n\in\mathbb{N}$. Then, using (\ref{m4}) with $v=\widetilde{u}$ we find that
\begin{equation*}
	\la A{u}_n,{u}_n-\widetilde{u}\ra
	\le \varphi(\widetilde{u})-\varphi(u_n)+\la f,{u}_n-\widetilde{u}\ra
\end{equation*}
or, equivalently,
\begin{equation}\label{z56}
\la A{u}_n-A\widetilde{u},{u}_n-\widetilde{u}\ra
	\le \la A\widetilde{u},\widetilde{u}-{u}_n\ra
	+\varphi(\widetilde{u})-\varphi(u_n)+\la  f,{u}_n-\widetilde{u}\ra .
\end{equation}

We use the strong monotonicity of the operator $A$, (\ref{A1}), to see that
\[
m_A\|u_n-\widetilde{u}\|_X^2\le
\la A\widetilde{u},\widetilde{u}-{u}_n\ra +\varphi(\widetilde{u})-\varphi(u_n)+\la  f,{u}_n-\widetilde{u}\ra .\]
Finally, we pass  to the upper limit in this inequality, use the weak convergence ${u}_n\rightharpoonup \widetilde{u}$ in $X$, and the weakly lower semicontinuity of the function $\varphi$ (guaranteed by assumption \eqref{fi}(a)) to infer that
\begin{equation}\label{z58}
	u_n\to \widetilde{u}\quad{\rm in}\quad X.	
\end{equation}
We now combine \eqref{z55} and \eqref{z58}
to see that $d(u_n,K)\to 0$ which concludes the proof of this step.

\medskip\noindent {\it Step
	{iii)  We  prove that
		any weakly convergent subsequence of the sequence $\{u_n\}$   converges to the solution $u$ of inequality $(\ref{1})$.}} This step is a direct consequence of Steps i), ii) and Theorem \ref{tm}.

\medskip\noindent {\it Step
	{iv)  We  prove that  the whole sequence $\{u_n\}$ converges to the solution $u$ of inequality $(\ref{1})$. }} To this end, we argue by contradiction and assume that  the convergence $u_n\to u$ in $X$ does not hold.  Then, there exists $\delta_0>0$ such that for all $k\in\mathbb{N}$ there exists $u_{n_k}\in X$ such that
\begin{equation}\label{m6}
	\|u_{n_k}-u\|_X\ge \delta_0.
\end{equation}
Note that the sequence $\{u_{n_k}\}$ is a subsequence of the sequence $\{u_n\}$ and, therefore, Step~i) and Lemma~\ref{lm} imply that it is  bounded in $X$. We now use a compactness argument to deduce that there exists a subsequence of the sequence $\{u_{n_k}\}$, again denoted by $\{u_{n_k}\}$, which is weakly convergent in $X$. Then, Step iii) guarantees that $u_{n_k}\to u$ as $k\to \infty$. We now pass at the limit when  $k\to \infty$ in \eqref{m6} and find $\delta_0\le 0$. This contradicts inequality $\delta_0>0$ and concludes the proof.
\hfill$\Box$

\section{A frictional contact problem}\label{s5}
\setcounter{equation}0

Theorem \ref{tp} can be used
in the study of various mathematical models which describe the equilibrium of elastic bodies in  frictional contact with an obstacle, the so-called foundation. It provides convergence results which lead to interesting mechanical interpretations. In this section
we introduce and study an example of such model and, to this end, we need some  notations and preliminaries.

Let $d\in\{2,3\}$. We denote by $\mathbb{S}^d$ the space of second order symmetric tensors on $\mathbb{R}^d$ and use  the notation $``\cdot"$, $\|\cdot\|$, $\bzero$ for the inner product, the norm and the zero element of the spaces
$\mathbb{R}^d$ and $\mathbb{S}^d$, respectively.
Let $\Omega\subset\mathbb{R}^d$ be a domain with a smooth boundary~$\Gamma$ divided into three
measurable disjoint parts $\Gamma_1$, $\Gamma_2$ and $\Gamma_3$ such that ${ meas}\,(\Gamma_1)>0$ and ${ meas}\,(\Gamma_3)>0$.
A generic point in $\Omega\cup\Gamma$ will be denoted by $\bx=(x_i)$.
We use the
standard notation for Sobolev and Lebesgue spaces associated to
$\Omega$ and $\Gamma$.
 In particular, we  use the spaces
\begin{eqnarray*}
	&&\hspace{-1.0cm}L^2(\Omega)^d=\Big\{\,\bv=(v_i)\ :\ v_i\in
	L^2(\Omega),\
	1\le i\le d\,\Big\},
	\\ [2mm]
	&&\hspace{-1.0cm}Q=L^2(\Omega)^{d\times
		d}_s=\Big\{\,\btau=(\tau_{ij}):\, \tau_{ij}=\tau_{ji}\in
	L^2(\Omega),
	\ 1\le i,\,j\le d\,\Big\}.
\end{eqnarray*}
These spaces  are real Hilbert spaces with the canonical inner products
\begin{eqnarray*}
	&&(\bu,\bv)_{L^2(\Omega)^d} = \int_\Omega u_i\,v_i\,dx=\int_\Omega \bu\cdot\bv\,dx\quad \forall\,
	\bu=(u_i),\,\bv=(v_i)\in L^2(\Omega)^d,\label{n1}\nonumber\\
	&&(\bsigma,\btau)_Q=\int_\Omega\sigma_{ij}\,\tau_{ij}\,dx=\int_\Omega\bsigma\cdot\btau\,dx
	\quad\forall\,\bsigma=(\sigma_{ij}),\,
	\btau=(\tau_{ij})
	\in Q\label{nq}\label{I2}
\end{eqnarray*}
and the associated norms  denoted by $\|\cdot\|_{L^2(\Omega)^d}$
and $\|\cdot\|_Q$, respectively. Here, the summation convention over repeated indices is used.

For an element $\bv\in H^1(\Omega)^d$ we still  write $\bv$ for the trace of
$\bv$ to $\Gamma$ and,
in addition, we consider the
space
\begin{eqnarray*}
	&&V=\{\,\bv=(v_i)\in H^1(\Omega)^d\ :\  \bv =\bzero\ \ {\rm on\ \ }\Gamma_1\,\},
\end{eqnarray*}
which is a real Hilbert space
endowed with the canonical inner product
\begin{equation*}
	(\bu,\bv)_V=
	(\bvarepsilon(\bu),\bvarepsilon(\bv))_Q
\end{equation*}
and the associated norm
$\|\cdot\|_V$. Here and below we use the symbol $\bvarepsilon(\bv)$
for the linear strain field, i.e.,
$\bvarepsilon(\bv)\in Q$ denotes the symmetric part of the gradient of $\bv\in V$.
We recall that, for an element $\bv\in V$,  the  normal and tangential components on $\Gamma$
are given by
$v_\nu=\bv\cdot\bnu$ and $\bv_\tau=\bv-v_\nu\bnu$, respectively.
We also  recall the trace inequality
\begin{equation}\label{trace}
	\|\bv\|_{L^2(\Gamma)^d}\leq d_0\|\bv\|_{V}\qquad \forall\,
	\bv\in V
\end{equation}
in which $d_0$ represents a positive constant.

For the inequality problem we consider in this section we use the data ${\cal F}$, $F_b$, $\fb_0$, $\fb_2$  which satisfy
the following conditions.
\begin{eqnarray}
	&&\label{Fc}
	\left\{ \begin{array}{ll} {\rm (a)\ } {\cal F}:\Omega\times
		\mathbb{S}^d\to \mathbb{S}^d.\\ [2mm]
		{\rm (b)\  There\ exists}\ L_{\cal F}>0\ {\rm such\ that}\\
		{}\qquad \|{\cal F}(\bx,\bvarepsilon_1)-{\cal F}(\bx,\bvarepsilon_2)\|
		\le L_{\cal F} \|\bvarepsilon_1-\bvarepsilon_2\|\\
		{}\qquad \quad\forall\,\bvarepsilon_1,\bvarepsilon_2
		\in \mathbb{S}^d,\ {\rm a.e.}\ \bx\in \Omega.\\ [2mm]
		{\rm (c)\  There\ exists}\ m_{\cal F}>0\ {\rm such\ that}\\
		{}\qquad ({\cal F}(\bx,\bvarepsilon_1)-{\cal F}(\bx,\bvarepsilon_2))
		\cdot(\bvarepsilon_1-\bvarepsilon_2)\ge m_{\cal F}\,
		\|\bvarepsilon_1-\bvarepsilon_2\|^2\\
		{}\qquad\quad \forall\,\bvarepsilon_1,
		\bvarepsilon_2 \in \mathbb{S}^d,\ {\rm a.e.}\ \bx\in \Omega.\\ [2mm]
		{\rm (d)\ The\ mapping}\ \bx\mapsto
		{\cal F}(\bx,\bvarepsilon)\ {\rm is\ measurable\ on\
		}\Omega,\\
		{}\qquad{\rm for\ any\ }\bvarepsilon\in \mathbb{S}^d.\\ [2mm]
		{\rm (e)\ The\ mapping\ } \bx\mapsto {\cal F}(\bx,\bzero)\ {\rm belongs\ to}\ Q.
	\end{array}\right.
	\label{IF}
\end{eqnarray}
\begin{eqnarray}
	&&\label{f0}\qquad F_b\in L^2(\Gamma_3),\qquad \fb_0\in L^2(\Omega)^d,\qquad\fb_2\in L^2(\Gamma_2)^d.\\ [3mm]
	&&\label{sm}\qquad d_0^2\,\|F_b\|_{L^2(\Gamma_3)} < m_{\cal F}.
\end{eqnarray}

\noindent
Recall that in (\ref{sm}) and below $d_0$ and $m_{\cal F}$ represent the constants introduced in \eqref{trace} and \eqref{Fc}, respectively. Moreover, we use $K$ for the set defined by
\begin{equation}
	\label{KK}K=\{\,\bv\in V\ :\ v_\nu \le 0\ \  \hbox{a.e. on}\
	\Gamma_3\,\}.
\end{equation}

Then, the inequality problem we consider  is the following.

\medskip\medskip\noindent
\medskip\noindent{\bf Problem}  ${\cal P}^c$. {\it Find  $\bu$
	such that}
\begin{eqnarray}\label{51}
	&&\bu\in K,\quad \int_{\Omega}{\cal F}\bvarepsilon(\bu)\cdot(\bvarepsilon(\bv)-\bvarepsilon(\bu))\,dx
       +\int_{\Gamma_3}F_b\,(\|\bv_\tau\|-\|\bu_\tau\|)\,da\\[2mm]
	&&\qquad\quad\ge
	\int_{\Omega}\fb_0\cdot(\bv-\bu)\,dx
	+\int_{\Gamma_2}\fb_2\cdot\gamma(\bv-\bu)\,da\quad\forall\,\bv\in K.\nonumber
\end{eqnarray}

\medskip

\noindent
Here, $\gamma\colon V\to L^2(\Gamma_3)^d$ is the trace operator.

Following the arguments in \cite{SMBOOK}, it can be shown that Problem $\cP^c$ represents the variational formulation of a mathematical model that describes the equilibrium of an elastic body  $\Omega$ which is acted upon by  external forces, is fixed on $\Gamma_1$, and is
in frictional contact on $\Gamma_3$ with a rigid foundation.
Here $\bu$ represents the displacement field, ${\cal F}$ is the elasticity operator,   $\fb_0$ and $\fb_2$ denote the density of  applied body forces and traction which act on the body and the surface $\Gamma_2$, respectively and $F_b$ is a given function, the friction bound.

Next, consider a sequence $\{\lambda_n\}\subset\real_+$ and a function $j_\nu$  which satisfy the following conditions.
\begin{eqnarray}
	&&\label{la}
\quad\ \lambda_n>0\quad\forall\, n\in\mathbb{N}, \ \lambda_n\to 0\ \ {\rm as}\ n\to\infty.	
\\ [3mm]
&&\label{jnu}
	\left\{
	\begin{array}{l}
		j_\nu \colon \Gamma_3 \times\mathbb{R} \to \mathbb{R}\ \mbox{is
			such that}
		\\[2mm]
		\ \ {\rm (a)}\ j_\nu(\cdot, r) \ \mbox{is measurable on} \ \Gamma_3 \ \mbox{for all} \ r \in \mathbb{R}
		\ \mbox{and there exists}\ \ \\
		\qquad \quad e_1 \in L^2(\Gamma_3) \ \mbox{such that} \
		j_\nu(\cdot, e_1(\cdot)) \in L^1(\Gamma_3).
		\\[2mm]
		\ \ {\rm (b)}\ j_\nu (\bx, \cdot) \ \mbox{is locally Lipschitz on} \ \real \ \mbox{for a.e.}\ \bx \in \Gamma_3.
		\\[2mm]
		\ \ {\rm (c)}\ | \partial j_\nu (\bx,r) | \le c_{0\nu} + c_{1\nu} | r | \ \mbox{for all}\ r \in \mathbb{R},
		\ \mbox{a.e.} \ \bx \in \Gamma_3 \\
		\qquad\quad \mbox{with}\ c_{0\nu}, c_{1\nu} \ge 0.
		\\[2mm]
		\ \ {\rm (d)}\ j_\nu^0(\bx, r; -r) \le d_\nu \left( 1 + |r| \right)
		\ \mbox{for all}\ r \in \mathbb{R}, \ \mbox{a.e.}\ \bx \in \Gamma_3\\
		\qquad\quad \mbox{with}\ d_\nu\ge 0.\\ [2mm]
		\ \ {\rm (e)} \ \ j_\nu(\bx,\cdot)\   \mbox{is regular 	for a.e.}\ \bx \in \Gamma_3.
	\end{array}
	\right.
\end{eqnarray}

\medskip

\noindent
Consider now  the following additional conditions.
\begin{eqnarray}
	&&\label{j1n}j^0_\nu(\bx,r;s-r)\le 0\qquad\forall\, r\in \real,\ s\le 0, \ \ \mbox{a.e.}\ \ \bx \in \Gamma_3. \\[2mm]
	&&\label{j2n}r\in \real,\quad j^0_\nu(\bx,r;s-r)\ge 0\quad\forall\, s\le 0,
	\ \mbox{a.e.}\ \ \bx \in \Gamma_3\ \Longrightarrow\ r\le 0. \\[2mm]
    &&\label{bdd} \forall\, s\in\real\ \  \exists\, \bv\in V \ \mbox{such\ that}\ v_\nu(\bx)=s \ \mbox{a.e.}\ \bx\in\Gamma_3.
\end{eqnarray}

\noindent Let us remark that if condition (\ref{j1n}) holds, then (\ref{jnu})(d) is also satisfied.
 Moreover, examples of sets $\Omega$, $\Gamma_1$ and $\Gamma_3$ for which condition (\ref{bdd}) is satisfied can be found in \cite{ST3}.

\smallskip

With these ingredients we
consider the following perturbation of Problem $\cP^c$.

\medskip\medskip\noindent{\bf Problem}  ${\cal P}^c_n$. {\it Find $\bu_n$ such that}				
\begin{eqnarray}\label{51n}
	&&\bu_n\in V,\quad \int_{\Omega}{\cal F}\bvarepsilon(\bu_n)\cdot(\bvarepsilon(\bv)-\bvarepsilon(\bu_n))\,dx\\[2mm]
	&&\qquad+\int_{\Gamma_3}F_b\,(\|{\bv}_\tau\|-\|{\bu_{n\tau}}\|)\,da+\frac{1}{\lambda_n}\int_{\Gamma_3}j_\nu^0(u_{n\nu};v_\nu-u_{n\nu})\,da\nonumber\\ [2mm]
	&&\qquad\qquad\ge
	\int_{\Omega}\fb_{0}\cdot(\bv-\bu_n)\,dx
	+\int_{\Gamma_2}\fb_{2}\cdot\gamma(\bv-\bu_n)\,da\quad\forall\,\bv\in V.\nonumber
\end{eqnarray}

The mechanical interpretation of Problem $\cP^c_n$ is similar to that of  Problem $\cP^c$. The difference arises in the fact that now the constraint $u_\nu\le 0$ is removed and, therefore, the contact is assumed to be with a deformable foundation. Here $j_\nu$  is a nonsmooth function which describes the reaction of the foundation, $\lambda_n$ is a deformability coefficient and $\frac{1}{\lambda_n}$ represents the stiffness of the foundation. For more comments see Remark \ref{R1} below.

\medskip

Our main result in this section,  is the following.

\begin{Theorem}\label{t5}
	Assume $(\ref{IF})$--$(\ref{sm})$, $(\ref{la})$--$(\ref{bdd})$. Then
	Problem $\cP^c$  has a unique solution  and, for each $n\in\mathbb{N}$, Problem $\cP^c_n$ has at least one solution.
	Moreover, if $\bu$ is the solution of Problem $\cP^c$ and
 $\{\bu_n\}$ is a~sequence of elements in $V$ such that, for each $n\in\mathbb{N}$, $\bu_n$ is a~solution of Problem $\cP_n^c$, then
 \begin{equation}\label{con}
 \bu_n\to \bu\qquad {\rm in}\ \ V\ \  {\rm as}\ \ n\to\infty.
 \end{equation}
\end{Theorem}

\medskip\noindent{\it Proof.}
	We start with some additional notation. First,  we
	consider the operator $A\colon V\to V^*$, the functions  $\varphi\colon V\to\R$, $j\colon L^2(\Gamma_3)^d\to\real$, and the element $\fb\in V^*$ defined as follows:
	\begin{eqnarray}
		&&\hspace{-14mm}
		\label{b1}\la A\bu,\bv\ra =\int_{\Omega}\cF\bvarepsilon(\bu)\cdot\bvarepsilon(\bv)\,dx,\\ [2mm]
		&&\hspace{-14mm}\label{b2}
		\varphi(\bv)=\int_{\Gamma_3} F_b\,\|{\bv}_\tau\|\,da,\\ [2mm]
		&&\hspace{-14mm}\label{b3}
		j(\bxi)=\int_{\Gamma_3} j_\nu(\xi_\nu)\,da,\\ [2mm]
		&&\hspace{-14mm}\label{b4}\la \fb,\bv\ra=\int_{\Omega}\fb_0\cdot\bv\,dx
		+\int_{\Gamma_2}\fb_2\cdot\gamma\bv\,da,
	\end{eqnarray}
	for all $\bu,\bv\in V$, $\bxi\in Y=L^2(\Gamma_3)^d$. Note that, since $j_\nu$ is regular, it follows that
	
\begin{equation}\label{2a}
j^0(\gamma\bu;\gamma\bv)=\int_{\Gamma_3}j_\nu^0(u_{\nu};v_\nu)\,da\qquad\forall\, \bu,\bv \in V.
\end{equation}

	Then it is easy to see that
	\begin{equation}\label{e1c}
		\left\{\begin{array}{l}
			\mbox{$\bu$ is a solution of Problem $\cP^c$ if and only if}\\ [3mm]
			\bu\in K, \quad \la A\bu,\bv-\bu\ra +\varphi(\bv)-\varphi(\bu)\ge \la\fb,\bv-\bu\ra\quad \forall\, \bv\in K.
		\end{array}\right.
	\end{equation}
	Moreover, for each $n\in\mathbb{N}$, the following equivalence holds: \begin{equation}\label{e2c}
		\left\{\begin{array}{l}
			\mbox{$\bu_n$ is a solution of Problem $\cP^c_n$ if and only if\qquad}\\ [3mm]
			\bu_n\in V, \quad \la A\bu_n,\bv-\bu_n\ra
			+\varphi(\bv)-\varphi(\bu_n)\\ [2mm]
			\qquad\qquad\qquad+\frac{1}{\lambda_n} j^0(\gamma\bu_n;\gamma\bv-\gamma\bu_n)\ge \la\fb,\bv-\bu_n\ra\quad \forall\, \bv\in V.
		\end{array}\right.
	\end{equation}

	\medskip\noindent
	Equivalence (\ref{e1c}) suggests us to use the abstract results in Sections \ref{s2} and \ref{s4} with $X=V$, $K$ defined by (\ref{KK}), $A$ defined by (\ref{b1}),
	$\varphi $ defined by (\ref{b2}) and $\fb$ given by (\ref{b4}). It is easy to see that in this case conditions $(\ref{K})$--$(\ref{f})$ are satisfied.
	For instance, using assumption (\ref{IF})  we see that
	\begin{eqnarray*}
		\la A\bu - A\bv,\bu -\bv\ra \geq m_{\cal F} \|\bu -\bv\|^{2}_{V},\qquad \|A\bu-A\bv\|_{V^*}\le L_{\cal F}\, {\|\bu-\bv\|_V}
	\end{eqnarray*}
	for all $\bu,\, \bv\in V$. Therefore, conditions (\ref{A}) and  (\ref{A1}) hold with $m_A=m_{\cal F}$ and $M_A=L_{\cal F}$, respectively.
	Conditions (\ref{fi}) and (\ref{fi1})  are also satisfied
	with $L_D= d_0^2 \|F_b\|_{L^2(\Gamma_3)}$, for any $D>0$.
	Therefore, we are in a position to apply Theorem \ref{t1}  in order to deduce the existence of a~unique solution of the variational inequality in (\ref{e1c}).
	Moreover, it is easy to see that the function $j$ defined by (\ref{b3}) inherits the properties of $j_\nu$ (\cite[Theorem 3.47]{MOSBOOK}). Hence, the solvability of the variational-hemivariational inequality in (\ref{e2c}) follows from Theorem~\ref{t2}.

Finally, we use assumption  \eqref{j1n},
\eqref{j2n} and  equality \eqref{2a} to see that conditions  \eqref{j1},
\eqref{j2} are satisfied.
Indeed, let $\bu\in V$ and $\bv\in K$.
Then, we have $v_\nu\le 0$ a.e. on $\Gamma_3$ and, therefore, assumption \eqref{j1n} implies that
$j_\nu^0(u_\nu; v_\nu-u_\nu)\le 0$ a.e. on $\Gamma_3$. We now use equality \eqref{2a}
to see that
\begin{equation}\label{r1}
	j^0(\gamma\bu;\gamma\bv-\gamma\bu)=\int_{\Gamma_3}j_\nu^0(u_\nu; v_\nu-u_\nu)\,da \le 0,
\end{equation}
which shows that condition (\ref{j1}) is satisfied.
Next, let $\bu\in V$, $\bv\in K$ and assume that
$j^0(\gamma\bu;\gamma\bv-\gamma\bu)\ge 0$ for all $\bv\in K$. We use \eqref{r1} to see that
\[j^0(\gamma\bu;\gamma\bv-\gamma\bu)=\int_{\Gamma_3}j_\nu^0(u_\nu; v_\nu-u_\nu)\,da=0. \]

Hence, because the integrand has nonpositive values,
we find that $j_\nu^0(u_\nu; v_\nu-u_\nu)=0$ a.e. on $\Gamma_3$. Since the previous inequality  holds for every $\bv\in K$, it follows  from assumptions (\ref{j2n}) and (\ref{bdd}) that $u_\nu\le 0$ a.e. on $\Gamma_3$ and, therefore, $\bu\in K$. This implies that
condition (\ref{j2}) is satisfied.

It follows from above that we are in a position to use Theorem \ref{tm} to deduce the convergence \eqref{con}.
These results combined with equivalences (\ref{e1c}) and (\ref{e2c}) allows us to conclude the proof of the theorem.
\hfill$\Box$

\begin{Remark}\label{R1}
Note that, among other ingredients, the variational inequality $(\ref{51})$ is obtained by using the Signorini contact condition in the form without gap, that is
\begin{equation*}\label{Sig}
	u_\nu\le 0,\quad \sigma_\nu\le 0,\quad \sigma_\nu u_\nu = 0\quad{\rm a.e.\ on}\ \ \Gamma_3.
\end{equation*}
Here and below $\sigma_\nu$ represents the normal stress on the contact boundary i.e., the normal component of $\bsigma\bnu$, where $\bsigma = {\cal F}\bvarepsilon(\bu)$. In contrast,  the variational-hemivariational inequality $(\ref{51n})$ is obtained by using the nonsmooth contact condition
\begin{equation}\label{nc}
	-\sigma_\nu\in\frac{1}{\lambda_n}\,\partial j_\nu(u_\nu)\quad{\rm a.e.\ on}\ \ \Gamma_3.
\end{equation}
This condition  represents a contact condition with normal compliance.
It describes the contact with a deformable foundation.
Here $\frac{1}{\lambda_n}$ can be interpreted as a stiffness coefficient of the foundation.  Indeed, in \eqref{nc} the penetration  is allowed but penalized.
\end{Remark}

\medskip
We now provide the following physical interpretation of Theorem \ref{t5}.
First, the existence and uniqueness part in the theorem proves the unique weak solvability of the frictional contact with a rigid foundation  and the weak solvability of the frictional contact with a deformable foundation.
Second,
the weak solution of the frictional contact problem with a rigid foundation   material  can be approached by the solution of the contact problem with a deformable foundation, with a large stiffness coefficient.

\medskip

We end this section with some examples of contact conditions which lead to
subdifferential conditions of the form (\ref{nc}).
Consider the normal compliance contact condition
\begin{equation}\label{G12.E0m}
	-\sigma_\nu=\frac{1}{\lambda_n}\, p_\nu(u_\nu)\quad{\rm on}\ \Gamma_3
\end{equation}
where $p_\nu \colon \real\to\real$ is a prescribed nonnegative
continuous function which vanishes when its argument is negative. Let $j_\nu \colon \real\to\real$ be the function defined by
\begin{eqnarray}
	&&\label{G.XJ1n}
	j_\nu(r) = \int_0^r p_\nu(s) \, ds \ \ \ \mbox{for all} \ \ r \in \real.
\end{eqnarray}

\noindent 
Then, we have
\begin{equation*}
	\displaystyle \partial j_\nu (r) = p_\nu (r) \ \ \ \mbox{for all} \ \ r \in \real
\end{equation*}
and,
therefore, it is easy to see that the contact condition
(\ref{G12.E0m}) is of the form (\ref{nc}).
Moreover, the conditions (\ref{j1n}) and (\ref{j2n}) reduce in this case to the following
\begin{eqnarray*}
	&&\label{j1'}p(r)\,(s-r)\le 0\qquad\forall\, r\in \real,\ s\leq 0,\\[2mm]
	&&\label{j2'} r\in \real,\quad p(r)\,(s-r)\ge 0\quad\forall\, s\leq 0\quad\Longrightarrow\ r\leq 0,
\end{eqnarray*}
respectively.

Concrete examples  lead to functions $j_\nu$ which satisfy
conditions (\ref{jnu}), (\ref{j1n}) and (\ref{j2n}) can be found in \cite{MOSBOOK}. Here we restrict ourselves to recall the following ones.

\begin{Example}\label{CC.Example1}
	Let $p_\nu \colon \real \to \real$ be the function given by
	\begin{equation*}\label{Ex1j}
		p_\nu(r)=a \, r_+ =  \left\{\begin{array}{ll} \   0\qquad &{\rm if}\ \ r<0,\\[2mm]
			\   ar  & {\rm if}\ \ r\ge 0,\end{array}\right.
	\end{equation*}
	with $a>0$. Then, using $(\ref{G.XJ1n})$ we have
	\begin{equation*}\label{Ex1g}
		j_\nu(r)=\left\{\begin{array}{ll} \ 0\qquad &{\rm if}\ \ r<0,\\[2mm]
			\ \displaystyle  \frac{a r^2}{2}  & {\rm if}\ \ r\ge
			0.\end{array}\right.
	\end{equation*}
	
	\noindent
	This contact condition corresponds to a
	linear dependence of the reactive force  with respect to the
	penetration and, therefore, it models a linearly elastic behaviour
	of the foundation.
\end{Example}

\begin{Example}\label{CC.Example2}
	Let  $p_\nu \colon \real \to \real$ be the function given by
	\begin{equation*}\label{Ex2j}
		p_\nu(r)=  \left\{\begin{array}{ll} \   0\qquad &{\rm if}\ \ r<0,\\[2mm]
			\   ar  & {\rm if}\ \ 0\le r\le l,\\[2mm]
			\   al  & {\rm if}\ \ r>l,
		\end{array}\right.
	\end{equation*}
	with $a>0$ and $l > 0$. Then, using $(\ref{G.XJ1n})$ we have
	
	\begin{equation*}\label{Ex2g}
		j_\nu(r)=  \left\{\begin{array}{ll} \   0\qquad\qquad &{\rm if}\ \ r<0,\\[2mm]
			\ \displaystyle \frac{a r^2}{2} & {\rm if}\ \ 0\le r\le l,\\[2mm]
			\ \displaystyle  a l r - \frac{a l^2}{2}\ \  & {\rm if}\ \ r>l.
		\end{array}\right.
	\end{equation*}
	

	\noindent	
	This contact condition corresponds to an
	elastic-perfect plastic behaviour of the foundation. The
	plasticity consists in the fact that when the penetration reach
	the limit $l$,  then the surface offers no additional resistance.
\end{Example}

\begin{Example}\label{CC.Example10}
	Let  $p_\nu \colon \real \to \real$ be the function
	given by
	\begin{equation*}\label{Ex6j}
		p_\nu(r)=  \left\{\begin{array}{ll} \ 0\qquad &{\rm if}\ \ r<0,\\[3mm]
			\   \displaystyle \frac{a+ e^{-b}}{b}\, r\quad  & {\rm if}\ \ 0\le r\le b,\\[5mm]
			\   \displaystyle e^{-r} + a   & {\rm if}\ \ r>b,
		\end{array}\right.
	\end{equation*}
	
	\noindent with $a \ge 0$, $b > 0$. Then, using $(\ref{G.XJ1n})$ we
	have
	
	\begin{equation*}\label{Ex6g}
		j_\nu(r)=  \left\{\begin{array}{ll} \   0\qquad\qquad &{\rm if}\ \ r<0,\\[3mm]
			\ \displaystyle
			\frac{a+e^{-b}}{2b} \, r^2 & {\rm if}\ \ 0\le r\le b,\\[4mm]
			\ \displaystyle a r - e^{-r} + \frac{(b+2) e^{-b} - ab}{2}\ \  &
			{\rm if}\ \ r>b.
		\end{array}\right.
	\end{equation*}

\noindent	
Note that in contrast to the
	previous two examples, here the function $p_\nu$ is not increasing
	and, therefore, the potential function $j_\nu$ is not a convex
	function.  This contact
	condition corresponds to an elastic-plastic behaviour of the
	foundation, with softening. The softening effect consists in the
	fact that, when the penetration reach the limit $b$, then the
	reactive force decreases.
\end{Example}

\section{Numerical Simulations} \label{s6}
In the final section of this work, we present the results of numerical simulations in the study of the frictional contact problem presented in Section \ref{s5}.
Our aim is to provide a numerical validation of the convergence result in Theorem  \ref{t5}.
To this end, we consider a sequence of quasistatic two-dimensional problems which based on Problem ${\cal P}^c_n$ where the parameter ${\lambda_n}$ is inversely proportional to the hardness of the foundation.
As $\lambda_n$ approaches zero, we expect the sequence of solutions converge to the solution of the Signorini Problem~${\cal P}^c$.

\subsection{Simulation Setup}

\begin{figure}
    \centering
    \includegraphics[width=0.5\textwidth]{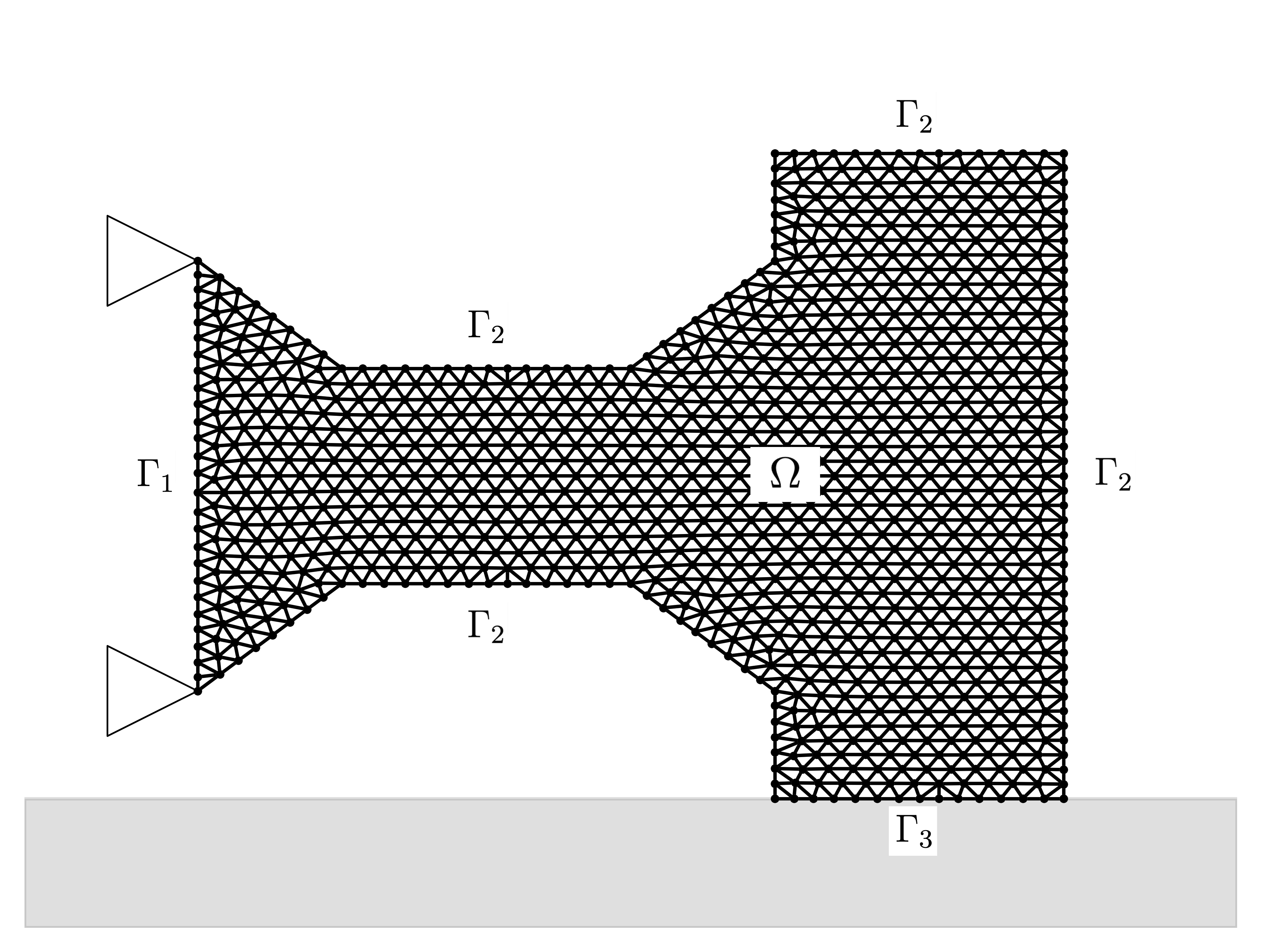}
    \caption{Mesh representation.}
\label{fig:mesh}
\end{figure}

The mesh corresponding to the physical setting we consider is depicted in Figure \ref{fig:mesh}.
For the simulations we use the finite element method with triangulation $\mathcal{T}$ of the domain $\Omega$, satisfying the necessary assumptions.
The adopted triangulation $\mathcal{T}$ is regular (according to \cite{ciarlet1978finite}) and shape-regular (according to \cite{braess2007finite}).
The spatial step for all simulations was chosen as $h = 1/32$.
The basis functions are well-known hat functions, i.e., piecewise linear polynomials ensuring continuity over the entire domain $\Omega$.
This discretization allows us to define the finite dimensional subspace $V^h \subset V$ and, furthermore, we introduce the discretized set $K^h = K \cap V^h$.

We now consider the following discrete version of Problem ${\cal P}^c$.

\medskip\noindent{\bf Problem}  ${\cal P}^{ch}$. {\it Find  $\bu^h$
	such that}
\begin{eqnarray*}\label{51h}
    &&\bu^h\in K^h,\quad \int_{\Omega}{\cal F}\bvarepsilon(\bu^h)\cdot(\bvarepsilon(\bv^h)-\bvarepsilon(\bu^h))\,dx+\int_{\Gamma_3}F_b\,(\|\bv^h_\tau\|-\|\bu^h_\tau\|)\,da\nonumber\\[2mm]
    &&\qquad\quad\ge
	\int_{\Omega}\fb_0\cdot(\bv^h-\bu^h)\,dx
	+\int_{\Gamma_2}\fb_2\cdot\gamma(\bv^h-\bu^h)\,da\quad\forall\,\bv^h\in K^h.\nonumber
\end{eqnarray*}
The discrete version of Problem ${\cal P}^c_n$ is as follows.

\medskip\medskip\noindent{\bf Problem}  ${\cal P}^{ch}_n$. {\it Find $\bu^h_n$ such that}				
\begin{eqnarray*}\label{51hn}
	&&\bu^h_n\in V^h,\quad \int_{\Omega}{\cal F}\bvarepsilon(\bu^h_n)\cdot(\bvarepsilon(\bv^h)-\bvarepsilon(\bu^h_n))\,dx\\[2mm]
	&&\qquad+\int_{\Gamma_3}F_b\,(\|{\bv^h}_\tau\|-\|{\bu^h_{n\tau}}\|)\,da+\frac{1}{\lambda_n}\int_{\Gamma_3}j_\nu^0(u^h_{n\nu};v^h_\nu-u^h_{n\nu})\,da\nonumber\\ [2mm]
	&&\qquad\qquad\ge
	\int_{\Omega}\fb_{0}\cdot(\bv^h-\bu^h_n)\,dx
	+\int_{\Gamma_2}\fb_{2}\cdot\gamma(\bv^h-\bu^h_n)\,da\quad\forall\,\bv^h\in V^h.\nonumber
\end{eqnarray*}

The mesh is illustrated in Figure~\ref{fig:mesh}, depicting a body fixed on the left side (the Dirichlet boundary condition) and placed on the foundation at the ``lowest" part where the contact boundary condition is defined. There is no gap between the body and the foundation. The influence of external forces $\boldsymbol{f}_2$ is omitted. For the sake of simplicity, we adopt a~dimensionless unit system. The material properties are characterized by the elasticity tensor with the components
\begin{eqnarray}
({\mathcal F}\boldsymbol{\omega})_{ij} = \frac{E\, \kappa}{(1 + \kappa)(1 - 2\kappa)}\,(\omega_{11} + \omega_{22})\,\delta_{ij} \, + \, \frac{E}{1 + \kappa}\,\omega_{ij} \label{elastic_tensor} \nonumber 
\end{eqnarray}
\noindent
for all $\boldsymbol{\omega}=(\omega_{ij}) \in \mathbb{S}^2$, $i,j = 1, 2$ with $E = 12 \cdot 10^3$ and $\kappa = 0.42$. We assume that the two-dimensional internal force field $\boldsymbol{f}_0$ is constant and equal to ($-0.2 \cdot 10^3$, $-0.8 \cdot 10^3$). The friction bound is constant, $F_b = 10$.
As a contact condition we utilized Example~\ref{CC.Example10} from the previous section with coefficients $a=b=0.1$.

\subsection{Simulation Results}
\begin{figure}[ht!]
    \vspace{13mm}
    \centering
    \begin{subfigure}[b]{0.49\textwidth}
        \centering
        \includegraphics[width=\textwidth]{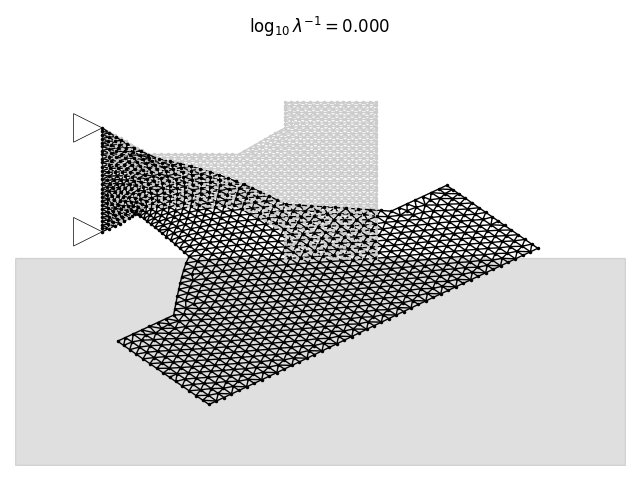}
        \caption{Almost no foundation.}
    \label{fig:results:0}
    \end{subfigure}
    \begin{subfigure}[b]{0.49\textwidth}
        \centering
        \includegraphics[width=\textwidth]{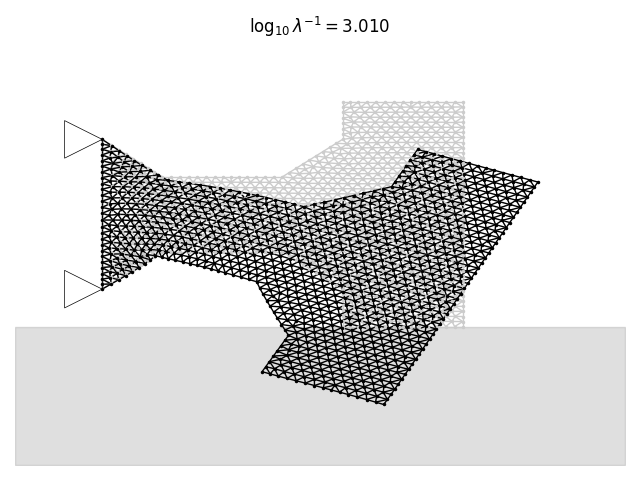}
        \caption{Soft foundation.}
    \label{fig:results:soft}
    \end{subfigure}
    \begin{subfigure}[b]{0.49\textwidth}
        \centering
        \includegraphics[width=\textwidth]{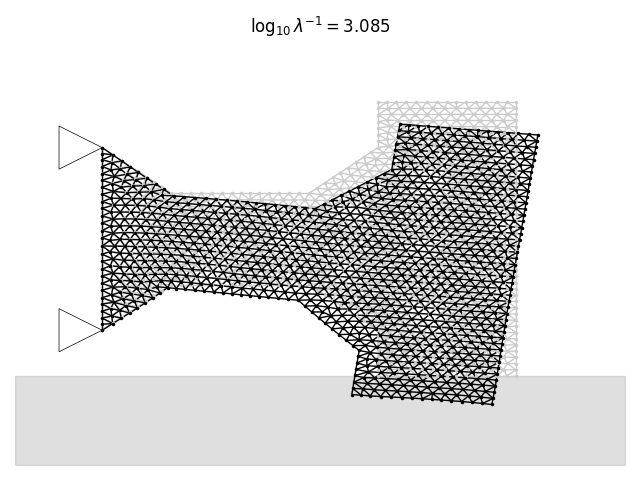}
        \caption{Moderately soft foundation.}
    \label{fig:results:m_soft}
    \end{subfigure}
    \begin{subfigure}[b]{0.49\textwidth}
        \centering
        \includegraphics[width=\textwidth]{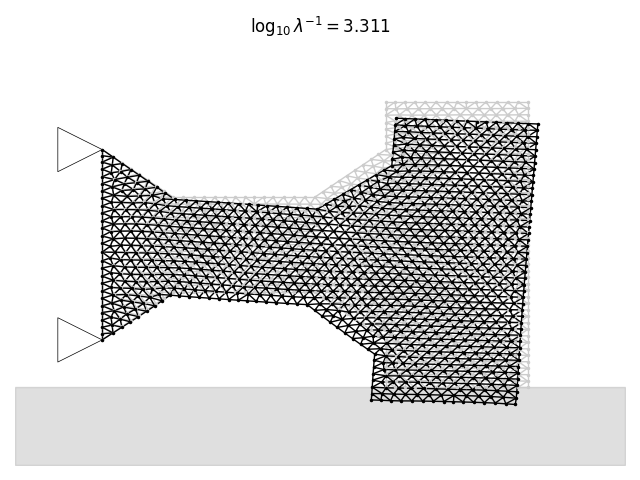}
        \caption{Moderately hard foundation.}
    \label{fig:results:m_hard}
    \end{subfigure}
    \begin{subfigure}[b]{0.49\textwidth}
        \centering
        \includegraphics[width=\textwidth]{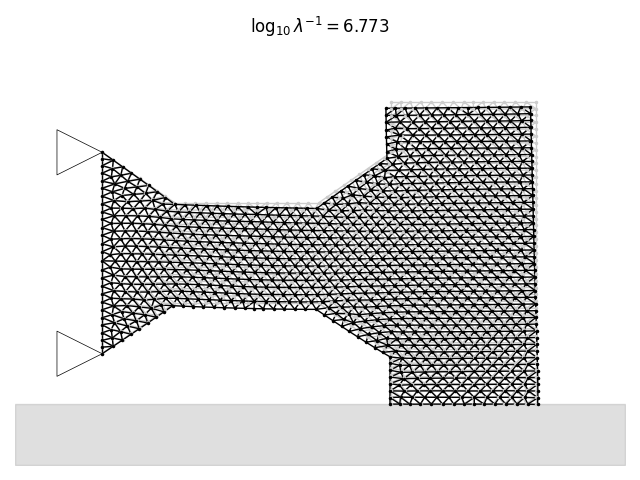}
        \caption{Hard foundation.}
    \label{fig:results:hard}
    \end{subfigure}
    \begin{subfigure}[b]{0.49\textwidth}
        \centering
        \includegraphics[width=\textwidth]{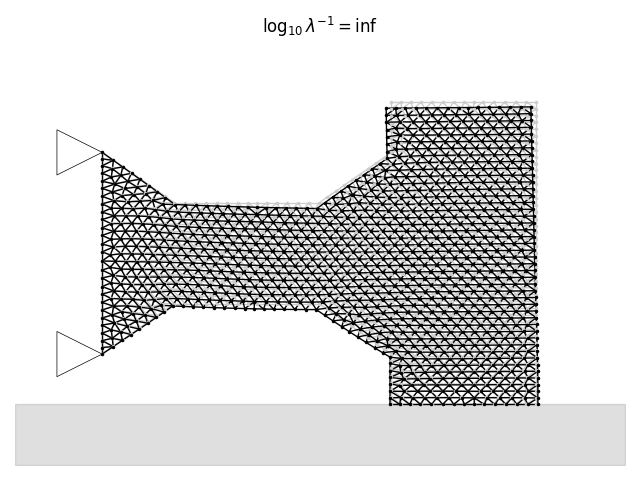}
        \caption{Rigid foundation.}
    \label{fig:results:rigid}
    \end{subfigure}
        \medskip
        \caption{Simulation results for selected $\lambda_n$ values.}
\label{fig:results}
\end{figure}

The results for selected values of $\lambda_n$ are presented in Figure~\ref{fig:results}. Above each figure, the approximate value of $\lambda_n$ is given. It should be noted that for readability, the values are expressed as the logarithm of the inversion of $\lambda_n$, meaning higher values indicate $\lambda_n$ itself approaching zero. In particular, the notation on the last figure specifies that it is a simulation in which a set of constraints $K$ was used instead of any value of $\lambda_n$. Figure~\ref{fig:results}(a) shows the body after applying force for almost zero foundation hardness, whereas Figure~\ref{fig:results}(b) illustrates the simulation result for a foundation with low hardness, allowing deep penetration. The second row displays the simulation result for a higher foundation hardness that shows resistance to the body.

\begin{figure}[ht]
    \centering
    \includegraphics[width=0.6\textwidth]{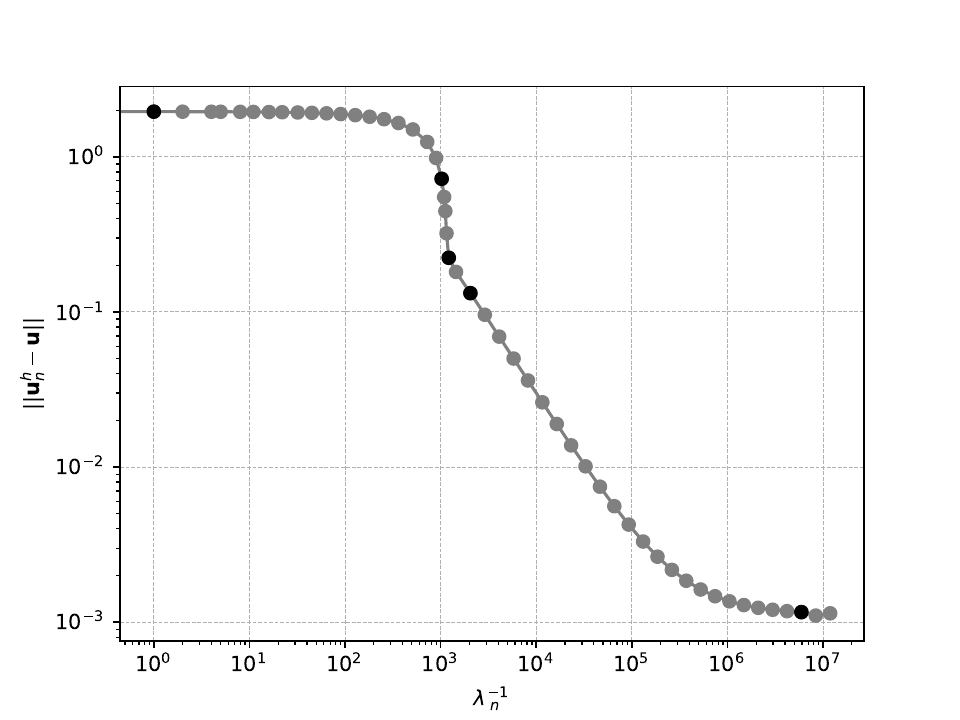}
    \caption{Convergence analysis showing the difference between solutions of Problem  ${\cal P}^{ch}_n$ for various $\lambda_n$ values and the reference solution. Highlighted points corresponds to plots in Figure~\ref{fig:results} (a)~-~(e).}
\label{fig:convergence}
\end{figure}

In Figure~\ref{fig:results}(e) the foundation is so hard that prevents downward movement, resembling the reference solution shown to the left. To analyze the convergence of solutions to the reference solution,  a special attention should be paid to Figure \ref{fig:convergence}, depicting the difference between the solution with the specified foundation hardness $\frac{1}{\lambda_n}$ ($x$-axis) and the reference solution.

The plot illustrates the displacement field difference compared to simulations with a~rigid foundation in terms of norm. It can be observed that as the value of $\lambda_n$ decreases, and consequently the stiffness of the foundation increases, the solutions converge. Particular attention should be paid to what happens when $\lambda_n$ takes values around $10^{-3}$. Note that, in our case of particular stiffness coefficients of the body as in the presented example, this is the point at which the integral over the contact boundary condition reaches a similar order of magnitude and begins to dominate in minimized functional with further growth. This implies that the solution to the series of Problems ${\cal P}^{ch}_n$ begins to significantly depend not only on the stiffness of the body itself but also on the resistance of the foundation. Hence, we can observe such a significant change in the course of solutions, followed by stabilization in convergence to the constrained solution.

The simulations are performed using the \textit{conmech} package \cite{OJB}. The package is a simulation tool developed in Python, designed specifically for addressing the complexities of contact mechanics problems across two-dimensional and three-dimensional bodies. It leverages the Finite Element Method (FEM) to facilitate the numerical resolution of static, quasistatic, and dynamic contact mechanics scenarios that defy analytical solutions. This package is adept to simulating a large range of physical phenomena including contact with various type of foundations, friction, adhesion and wear, among others. Built to be almost entirely self-sufficient, \textit{conmech} only requires basic Python libraries. Its modular design not only ensures that it covers a wide array of contact mechanics applications but also allows for future expansions to incorporate new physical models, making it an invaluable tool for researchers and engineers in the field. The package is available under the GPL-3.0 license, promoting open-source collaboration and innovation. To solve the sequence of unconstrained variational-hemivariational inequalities, an optimization approach utilizing the Powell method is employed \cite{JOB}. Conversely, Problem ${\cal P}^{ch}$ is addressed using optimization methods with constraints.

In conclusion, the presented simulations successfully validate the theoretical results for the frictional contact phenomenon. The systematic analysis of varying base hardness provides insights into the behavior of the body under external forces. Further exploration and refinement of the model can contribute to advancements in understanding and predicting frictional contact scenarios.

\section*{Acknowledgement}

\indent
This project has received funding from the European Union’s Horizon 2020
Research and Innovation Programme under the Marie Sklodowska-Curie
Grant Agreement No 823731 CONMECH.
The first two authors are supported by
National Science Center, Poland, under project OPUS no. 2021/41/B/ST1/01636.
The first author of the publication received an incentive scholarship from the funds of the program Excellence Initiative - Research University at the Jagiellonian University in Krakow.


\begin{thebibliography}{99}
	
	
	\bibitem{BC}
	C. Baiocchi and A. Capelo, {\em Variational and Quasivariational
		Inequalities: Applications to Free-Boundary Problems}, John Wiley,
	Chichester, 1984.
	
	
	
	\bibitem{braess2007finite}
    D. Braess,  \textit{Finite Elements – Theory, Fast Solvers, and Applications in Solid Mechanics}, 3rd edn. Cambridge University Press, New York, 2007.
	
	
	\bibitem{B}
	H. Br\'ezis, Equations et in\'equations non lin\'eaires dans les
	espaces vectoriels en dualit\'e, {\em Ann. Inst. Fourier Grenoble} {\bf 18} (1968), 115--175.
	
	
	
	\bibitem{C} A. Capatina, {\it Variational Inequalities
		Frictional Contact Problems}, Advances in Mechanics and Mathematics, Vol. 31, Springer, New York, 2014.
	

        \bibitem{ciarlet1978finite}
        P.G. Ciarlet, \textit{The Finite Element Method for Elliptic Problems}. North-Holland, Amsterdam, 1978.
	
\bibitem{Cl}	
	 F. H. Clarke,
    {\em Optimization and Nonsmooth Analysis\/},
     Wiley Interscience, New York, 1983.
	
	\bibitem{DL}
	G. Duvaut and J.-L. Lions, {\em Inequalities in Mechanics and
		Physics\/}, Springer-Verlag, Berlin, 1976.
	
	\bibitem{EJK}
	C. Eck, J. Jaru\v sek and M. Krbec, {\em Unilateral  Contact
		Problems: Variational Methods and Existence Theorems}, Pure and Applied Mathematics {\bf 270}, Chapman/CRC Press, New York, 2005.
	
	
	
	\bibitem{GMOT} C. Gariboldi, S. Migorski, A. Ochal and D. A. Tarzia, Existence, comparison, and convergence results for a class of elliptic hemivariational inequalities, {\it Applied Mathematics \& Optimization} {\bf 84} (Suppl 2) (2021), S1453--1475.
	
	\bibitem{GOST}
	C. Gariboldi, A. Ochal, M. Sofonea and D. Tarzia, Modelling, A Convergence Criterion for Elliptic Variational Inequalities,
{\it Applicable Analysis} {\bf 103(10)} (2024), 1810-1830
	
	
	
	
	
	\bibitem{G}
	R. Glowinski, {\em Numerical Methods for Nonlinear Variational
		Problems\/}, Springer-Verlag, New York, 1984.




\bibitem{Han}
W. Han,
A Revisit of Elliptic Variational-Hemivariational Inequalities,
{\it Numerical Functional Analysis and Optimization} {\bf 42(4)} (2021), 371-–395.



\bibitem{HS-Acta} W. Han and M. Sofonea, Numerical analysis of hemivariational inequalities
in Contact Mechanics, {\it Acta Numer.} (2019), 175--286.	
	
	\bibitem{HS}
	W. Han and M. Sofonea, {\it Quasistatic  Contact Problems in
		Viscoelasticity and Viscoplasticity}, Studies in Advanced
	Mathematics {\bf 30}, American Mathematical Society, Providence,
	RI--International Press, Somerville, MA, 2002.
	
	
	
	\bibitem{HHNL}
	I. Hlav\'a\v{c}ek, J. Haslinger, J. Nec\v{a}s and J. Lov\'{\i}\v{s}ek,
	{\em Solution of Variational Inequalities in Mechanics},
	Springer-Verlag, New York, 1988.

        \bibitem{JOB}
        M. Jureczka, A. Ochal and P. Bartman, {\it A nonsmooth optimization approach for time-dependent hemivariational inequalities}, Nonlinear Analysis: Real World Applications {\bf 73} (2023), 103871 .  
	
	
	
	\bibitem{KO} N. Kikuchi and J.T. Oden, {\em Contact Problems in
		Elasticity: A Study of Variational Inequalities and Finite Element
		Methods}, SIAM, Philadelphia, 1988.
	
	
	\bibitem{Kind-St}
	D. Kinderlehrer and G. Stampacchia, {\em An Introduction to
		Variational Inequalities and their Applications}, Classics in
	Applied Mathematics {\bf 31}, SIAM, Philadelphia, 2000.
	


\bibitem{LM}
Z.H. Liu, D. Motreanu,
A class of variational–hemivariational inequalities of elliptic type,
{\it Nonlinearity} {\bf 23(7)} (2010), 1741.



\bibitem{MNZ}
S. Migórski, V.T. Nguyen, S.D. Zeng,
Nonlocal elliptic variational-hemivariational inequalities,
{\it J. Integral Equations Applications} {\bf 32(1)} (2020), 51--58.


\bibitem{MYZ}
S. Migórski, J.-C. Yao, S.D. Zeng,
A class of elliptic quasi-variational–hemivariational inequalities with applications,
{\it Journal of Computational and Applied Mathematics} {\bf 421(15)} (2022), 114871.


\bibitem{MOSBOOK}
S. Mig\'orski, A. Ochal and M. Sofonea,
{\it Nonlinear Inclusions and Hemivariational Inequalities.
	Models and Analysis of Contact Problems},
Advances in Mechanics and Mathematics \textbf{26},
Springer, New York, 2013.

	
\bibitem{MOS30}
S. Mig\'orski, A. Ochal and M. Sofonea,
A class of variational-hemi\-va\-ria\-tio\-nal in\-equa\-li\-ties in reflexive Banach spaces,
{\it  J. Elasticity} {\bf 127} (2017), 151--178.




\bibitem{NP}
Z. Naniewicz and P. D. Panagiotopoulos,
{\it Mathematical Theory of Hemivariational Inequalities
	and Applications}, Marcel Dekker, Inc., New York, Basel,
Hong Kong, 1995.

\bibitem{OJB} A. Ochal, M. Jureczka and P. Bartman, A survey of numerical methods for hemivariational inequalities with applications to Contact Mechanics, {\it Communications in Nonlinear Science and Numerical Simulation} {\bf 114} (2022), 106563.

\bibitem{Pana1}
P.D. Panagiotopoulos, Nonconvex problems of semipermeable media
and related topics, {\em Z. Angew. Math. Mech. (ZAMM)} {\bf 65}
(1985), 29--36.

\bibitem{Pa} P.D. Panagiotopoulos,
{\em Inequality Problems in Mechanics and Applications\/},
Birkh\"{a}user, Boston, 1985.

\bibitem{P}
P. D. Panagiotopoulos, {\it Hemivariational Inequalities,
	Applications in Mechanics and Engineering}, Springer-Verlag,
Berlin, 1993.



\bibitem{PK} Z. Peng and  K. Kunisch, Optimal control of elliptic variational-hemivariational inequalities,
{\it J. of Optim. Theory Appl.} {\bf 178} (2018), 1--25.



\bibitem{S} M. Sofonea, {\it Well-posed Nonlinear Problems. A Study of Mathematical Models of Contact}, Advances in Mechanics and Mathematics {\bf 50}, Birkh\"auser, Cham, 2023.
	
	
\bibitem{SMBOOK} M. Sofonea and S. Mig\'orski,
	{\it Variational-Hemivariational Inequalities with Applications}, Pure and Applied Mathematics, Chapman \& Hall/CRC
	Press, Boca Raton-London, 2018.	

\bibitem{ST3}  M. Sofonea and D. A. Tarzia, Tykhonov well-posedness of a heat transfer problem with unilateral Constraints, {\it Applications of Mathematics} {\bf 67(2)} (2022), 167--197.
	
	
\end{thebibliography}
\end{document}